\documentclass[11pt]{amsart}
\usepackage{amssymb,amsmath,epsfig,graphics}
\theoremstyle{plain}
\newtheorem{thrm}{Theorem}[section]
\newtheorem{lemma}[thrm]{Lemma}
\newtheorem{prop}[thrm]{Proposition}
\newtheorem{cor}[thrm]{Corollary}
\newtheorem{rmrk}[thrm]{Remark}
\newtheorem{dfn}[thrm]{Definition}

\numberwithin{equation}{section} \numberwithin{figure}{section}

\begin{document}
\newcommand{\SL}{\mathcal L^{1,p}(\Om)}
\newcommand{\Lp}{L^p(\Omega)}
\newcommand{\CO}{C^\infty_0(\Omega)}
\newcommand{\Rn}{\mathbb R^n}
\newcommand{\Rm}{\mathbb R^m}
\newcommand{\R}{\mathbb R}
\newcommand{\Om}{\Omega}
\newcommand{\Hn}{\mathbb H^n}
\newcommand{\HH}{\mathbb H^1}
\newcommand{\eps}{\epsilon}
\newcommand{\BVX}{BV_H(\Omega)}
\newcommand{\IO}{\int_\Omega}
\newcommand{\bG}{\boldsymbol{G}}
\newcommand{\bg}{\mathfrak g}
\newcommand{\p}{\partial}
\newcommand{\Xnu}{\overset{\rightarrow}{ X_\nu}}
\newcommand{\nuX}{\boldsymbol{\nu}^H}
\newcommand{\Up}{\bN^H}
\newcommand{\n}{\boldsymbol \nu}
\newcommand{\sigmau}{\boldsymbol{\sigma}^u_H}
\newcommand{\di}{\nabla_{i}^{H,\mS}}
\newcommand{\one}{\nabla_{1}^{H,\mS}}
\newcommand{\two}{\nabla_{2}^{H,\mS}}
\newcommand{\del}{\nabla^{H,\mS}}
\newcommand{\delXY}{\nabla^{H,\mS}_X Y}
\newcommand{\nui}{\nu^H_i}
\newcommand{\nuj}{\nu^H_j}
\newcommand{\dej}{\nabla_{j}^{H\mS}}
\newcommand{\cx}{\boldsymbol{c}^{H,\mathcal S}}
\newcommand{\sx}{\sigma_H}
\newcommand{\lx}{\mathcal L_H}
\newcommand{\pb}{\overline p}
\newcommand{\qb}{\overline q}
\newcommand{\ob}{\overline \omega}
\newcommand{\sbar}{\overline s}
\newcommand{\nuu}{\boldsymbol \nu_{H,u}}
\newcommand{\nuv}{\boldsymbol \nu_{H,v}}
\newcommand{\Bl}{\Bigl|_{\lambda = 0}}
\newcommand{\Bs}{\Bigl|_{s = 0}}
\newcommand{\mS}{\mathcal S}
\newcommand{\delh}{\Delta_H}
\newcommand{\delinf}{\Delta_{H,\infty}}
\newcommand{\nabh}{\nabla^H}
\newcommand{\delp}{\Delta_{H,p}}
\newcommand{\mO}{\mathcal O}
\newcommand{\delhs}{\Delta_{H,\mS}}
\newcommand{\dla}{\Delta_{H,\mS^\la}}
\newcommand{\lhs}{\hat{\Delta}_{H,\mS}}
\newcommand{\bN}{\boldsymbol{N}}
\newcommand{\bnu}{\boldsymbol \nu}
\newcommand{\la}{\lambda}
\newcommand{\nup}{(\nuX)^\perp}
\newcommand{\nuup}{\boldsymbol \nu^\perp_{H,u}}
\newcommand{\nuvp}{\boldsymbol \nu^\perp_{H,v}}
\newcommand{\ep}{\epsilon}
\newcommand{\e}{\mathfrak e}
\newcommand{\E}{\mathfrak E}
\newcommand{\BH}{\Bigl|_{H_g}}
\newcommand{\om}{\omega}
\newcommand{\se}{\sqrt{\epsilon}}
\newcommand{\oX}{\overline X}
\newcommand{\oY}{\overline Y}
\newcommand{\ou}{\overline u}
\newcommand{\fv}{\mathcal V^{H}_I(\mS;\mathcal X)}
\newcommand{\sv}{\mathcal V^{H}_{II}(\mS;\mathcal X)}
\newcommand{\nh}{\nabla^H}
\newcommand{\sul}{\Delta_H}
\newcommand{\sij}{\sum_{i,j=1}^m}


\title[Sub-Riemannian calculus and monotonicity of the perimeter, etc.]
{Sub-Riemannian calculus and monotonicity of the perimeter for
graphical strips}

\author{D. Danielli}
\address{Department of Mathematics\\Purdue University \\
West Lafayette, IN 47907} \email[Donatella
Danielli]{danielli@math.purdue.edu}
\thanks{First author supported in part by NSF grant CAREER DMS-0239771}

\author{N. Garofalo}
\address{Department of Mathematics\\Purdue University \\
West Lafayette, IN 47907} \email[Nicola
Garofalo]{garofalo@math.purdue.edu}
\thanks{Second author supported in part by NSF Grant DMS-0701001}

\author{D. M. Nhieu}
\address{Department of Mathematics \\
San Diego Christian College \\
2100 Greenfield dr\\
El Cajon CA 92019} \email[Duy-Minh Nhieu]{dnhieu@sdcc.edu}

%
%
\date{\today}

\maketitle

\baselineskip 16pt



\section{\textbf{Introduction}}

In recent years the study of surfaces of constant horizontal mean
curvature $\mathcal H$ (to be defined below) in sub-Riemannian
spaces has seen an explosion of interest. Similarly to the classical
situation, this interest has provided a strong stimulus for the
development of a corresponding geometric measure theory. For a
partial account of such surge of activity the reader should consult
 \cite{Pa1}, \cite{Pa2}, \cite{CDG},
\cite{KR}, \cite{E1}, \cite{E2}, \cite{E3}, \cite{Gro}, \cite{GN},
\cite{Be}, \cite{DS}, \cite{DGN1}, \cite{AK1}, \cite{AK2},
\cite{CS1}, \cite{A}, \cite{FSS1}, \cite{Ma1}, \cite{FSS2},
\cite{Ma2}, \cite{CMS}, \cite{FSS3}, \cite{BRS}, \cite{DGN4},
\cite{DGN5}, \cite{DGN5}, \cite{LR}, \cite{LM}, \cite{FSS4}, \cite{Ma3},
\cite{CS2} \cite{P1}, \cite{P2}, \cite{GP1}, \cite{CG}, \cite{CHMY},
\cite{CH}, \cite{HP}, \cite{HP2}, \cite{RR}, \cite{BC}, \cite{Se},
\cite{Se2}, \cite{Mo}.

In this context, the Heisenberg group $\Hn$ occupies a central
position, especially in connection with the sub-Riemannian Bernstein
and isoperimetric problems. We recall that $\Hn$ is the stratified
nilpotent Lie group whose (real) underlying manifold is $\R^{2n+1}$
with the non-Abelian group law inherited by the complex product in
$\mathbb C^{n+1}$
\[
(x,y,t)\cdot (x',y',t')\ =\ \left(x+x',y+y',t+t' + \frac{1}{2}
(<x,y'> - <x',y>)\right).
\]

If we set $p = (x,y,t), p'=(x',y',t')\in \R^{2n+1}$, define the
left-translation map by $L_p(p') = p\circ p'$, and we indicate with
$L_p^*$ its differential, then the Lie algebra of all left-invariant
vector fields in $\Hn$ is spanned by the $2n+1$ vector fields
\[
X_i = L_p^*(\p_{x_i}) = \p_{x_i} - \frac{y_i}{2} \p_t,\ \ X_{n+i} =
L_p^*(\p_{y_i}) = \p_{y_i} + \frac{x_i}{2} \p_t,\ \ T =
L_p^*(\p_{t}) = \p_t,
\]
where $i=1,...,n$. We note the important commutation relations
$[X_i,X_{n+j}] = T$, $i,j=1,...,n$. They guarantee that the vector
fields $X_1,...,X_{2n}$ suffice to generate the whole Lie algebra,
and therefore the Heisenberg group is a stratified nilpotent Lie
group of step two, see \cite{Fo}, \cite{S}, \cite{BLU}. Such group
is in fact the basic model of such sub-Riemannian manifolds, and it
plays in this context much the same role played by $\Rn$ in
Riemannian geometry. The first Heisenberg group $\HH$ is obtained
when $n=1$. If we indicate with $p = (x,y,t)\in \R^3$ a generic
point of its underlying manifold, then the generators of its (real)
Lie algebra are the two vector fields
\[
X_1 = L_p^*(\p_{x})= \p_x - \frac{y}{2} \p _t,\ \ \ X_2 =
L_p^*(\p_{y}) =\p_y + \frac{x}{2} \p _t,
\]
and we clearly have $[X_1,X_2] = T = L_p^*(\p_{t})= \p_t$.

To introduce the results in this paper we recall that one of the
most fundamental properties of classical minimal surfaces $\mathcal
S\subset \Rm$ is the following well-known monotonicity theorem, see
\cite{MS}, and also \cite{Si}, \cite{MM}, \cite{CM}

\begin{thrm}\label{T:monoeuclid}
Let $\mathcal S\subset \R^m$ be a $C^2$ hypersurface, with $H$ being
its mean curvature, then for every fixed $p\in \mathcal S$ the
function
\begin{equation}\label{vg}
r\ \to\ \frac{H_{m-1}(\mathcal S \cap B_e(p,r))}{r^{m-1}}\ +\
\int_0^r \frac{m-1}{t^{m-1}}\ \int_{\mathcal S \cap B_e(p,t)} |H|\
dH_{m-1}\ dt\ ,
\end{equation}
is non-decreasing. In particular, if $\mS$ is minimal, i.e., if
$H\equiv 0$, then
\begin{equation}\label{vg2}
r\ \to\ \frac{H_{m-1}(\mathcal S \cap B_e(p,r))}{r^{m-1}}
\end{equation}
is non-decreasing.
\end{thrm}

In \eqref{vg}, \eqref{vg2} we have denoted by $H_{m-1}$ the
$(m-1)$-dimensional Hausdorff measure in $\Rm$. Theorem
\ref{T:monoeuclid} has many deep implications. It says, in
particular, that minimal hypersurfaces have maximum volume growth at
infinity, i.e., there exists $c_m>0$ such that $H_{m-1}(S\cap
B(p,r)) \geq\ c_m r^{m-1}$ as $r\to \infty$.

In this paper we are interested in related growth properties of the
sub-Riemannian volume on a $H$-minimal surface in $\HH$.  By
$H$-minimal we mean a $C^2$ oriented hypersurface $\mS\subset \Hn$
such that its horizontal mean curvature $\mathcal H$ vanishes
identically on $\mS$. The sub-Riemannian volume instead is the
so-called horizontal perimeter, see Section \ref{S:IBP} for its
definition and main properties. We should say right upfront that,
despite the efforts of several workers, the monotonic character of
the sub-Riemannian volume continues to represent a fundamental open
question.

The main obstacle so far has been represented by finding an
appropriate substitute of some basic properties such as, for
instance, the following elementary, yet fundamental fact from
Riemannian geometry. Consider in $\Rm$ the radial vector field
$\zeta(x) = \sum_{i=1}^m x_i \p_{x_i}$, then on \emph{any} $C^2$
hypersurface $\mS\subset \Rm$, one has
\begin{equation}\label{fi} div_\mS \zeta\ \equiv\ m-1,
\end{equation} where we have indicated with $div_\mS$ the Riemannian
divergence on $\mS$. The elementary identity \eqref{fi} has many
deep implications, and one could safely claim that behind most
fundamental results from the classical theory of minimal surfaces
there is \eqref{fi}. For instance, Theorem \ref{T:monoeuclid} and
the Sobolev inequalities on minimal surfaces \cite{MS} are
consequences (highly non-trivial, of course) of \eqref{fi}. The
number $m-1$ in the right-hand side of \eqref{fi} is dimensionally
correct since the standard volume form $\sigma$ on a hypersurface in
$\R^m$ scales according to the rule
\[
\sigma(\delta_\lambda(E)) = \lambda^{m-1} \sigma(E),\ \ \ E\subset
\mS,
\]
where $\delta_\lambda(x) = \lambda x$ represent the isotropic
dilations in $\R^m$.

In sub-Riemannian geometry, however, the correct dimension is
dictated by the non-isotropic dilations of the ambient non-Abelian
group, and this seemingly natural fact becomes a source of great
complications. For instance, given a $C^1$ hypersurface $\mS\subset
\Hn$, and indicating with $\sigma_H$ the horizontal perimeter on
$\mS$ (for its definition we refer the reader to Section
\ref{S:IBP}), then one has

\begin{equation}\label{scalings}
\sigma_H(\delta_\lambda(E))\ =\ \lambda^{Q-1} \sigma_H(E),\ \ \
E\subset \mS,
\end{equation}
where $\delta_\lambda(x,y,t) = (\lambda x,\lambda y, \lambda^2 t)$
indicates the non-isotropic dilations in $\Hn$ associated with the
grading of its Lie algebra. Here, the number $Q = 2n+2$ represents
the homogeneous dimension of $\Hn$ associated with the dilations
$\{\delta_\lambda\}_{\lambda>0}$. Thus for instance, when $n=1$, we
have $Q = 4$.

Guided by the analogy with \eqref{fi} one would like to find a
horizontal vector field $\zeta$ in $\Hn$ whose sub-Riemannian
divergence on $\mS$ (to be precisely defined below) satisfy the
equation
\begin{equation}\label{wrong}
div_{H,\mS} \zeta = Q -1.
\end{equation}

Such attempt would not possibly work however, for several reasons
which are all connected to one another. First of all, the
integration by parts formula in which one would like to use such a
$\zeta$ contains a corrective term which is produced by the above
mentioned non trivial commutation relations which connect the
generators of the Lie algebra of $\Hn$. Secondly, one should not
forget that not only the radial vector field $\zeta$ satisfies
\eqref{fi}, but it also possess the equally important property that
\begin{equation}\label{growth}
\underset{x\in\mS\cap B(0,r)}{\sup}\ |<\zeta(x) ,\nabla^\mS |x|>|\
\leq\ r,
\end{equation}
where $\nabla^\mS$ indicates the Riemannian gradient on $\mS$.
Because of these obstructions, there has been no progress so far on
the question of the monotonic character of sub-Riemannian minimal
surfaces.

One of the main contributions of the present paper is a monotonicity
formula for an interesting class of $H$-minimal surfaces in $\HH$,
the so-called \emph{graphical strips}. Such surfaces were introduced
in the work \cite{DGNP}, where they played a crucial role in the
solution of the sub-Riemannian Bernstein problem in $\HH$. Our main
result hinges on the discovery that, despite the original evidence
against it, for such class of surfaces the generator of the
non-isotropic group dilations in $\HH$ provides a valid replacement
of the radial vector field in $\R^m$. This sentence must, however,
be suitably interpreted, in the sense that things do not work so
simply. What we mean by this is that the horizontal integration by
parts formulas from \cite{DGN4} (see also \cite{DGN2}) which
constitute the sub-Riemannian counterpart of the classical
integration by parts formulas on hypersurfaces (for these, see e.g.
\cite{MM}, \cite{Si}, \cite{CM}), do not suffice. They need to be appropriately
intertwined with a twisted vertical integration by parts formula
also discovered in \cite{DGN4}. Both such formulas have played a
pervasive role in the establishment of a general second variation
formula for the horizontal perimeter. To state our main result we
recall the relevant definition.

\begin{dfn}\label{D:gs}
We say that a $C^1$ surface $\mS\subset \HH$ is a \emph{graphical
strip} if there exist an interval $I\subset \R$, and $G \in C^1(I)$,
with $G'\geq 0$ on $I$, such that, after possibly a left-translation
and a rotation about the $t$-axis, then either
\begin{equation}\label{ceI}
\mathcal S\ =\ \{(x,y,t)\in  \HH \mid  (y,t) \in \R \times I , x = y
G(t)\}\ ,
\end{equation}
or
\begin{equation}\label{ceII}
\mathcal S\ =\ \{(x,y,t)\in  \HH \mid  (x,t) \in \R \times I , y = -
x G(t)\}\ .
\end{equation}
If there exists $J\subset I$ such that $G'>0$ on $J$, then we call
$\mS$ a \emph{strict graphical strip}.
\end{dfn}

When the interval $I$ can be taken to be the whole real line, then
we call $\mS$ an entire graphical strip (strict, if $G'>0$ on some
$J\subset \R$).

The main result of this paper is the following theorem.

\begin{thrm}\label{T:monogsintro}
Let $\mS\subset \HH$ be a $C^2$ graphical strip, and denote by
$\sigma_H$ the sub-Riemannian volume form, or horizontal perimeter,
on $\mS$. For every $p_0=(0,0,t_0)\in \mS$ the function
\[
r\ \to\ \frac{\sigma_H(\mS\cap B(p_0,r))}{r^{Q-1}},\ \ \ \ r>0,
\]
is monotone non-decreasing. Moreover, there exists $\omega>0$ such
that
\[
\sigma_H(\mS\cap B(p_0,r)) \geq \omega r^{Q-1},\ \ \ \text{for
every}\ r>0.
\]
\end{thrm}

In the statement of Theorem \ref{T:monogsintro} we have denoted by
$B(p_0,r) = \{p\in \Hn \mid d(p,p_0)<r\}$, where $d(p,p_0) =
N(p_0^{-1} p)$ represents the gauge distance on $\Hn$ defined via
the Koranyi-Folland gauge function $N(p) = (|z|^4 + 16 t^2)^{1/4}$,
$p = (z,t)\in \Hn$.

The proof of Theorem \ref{T:monogsintro} is inspired to the ideas
set forth in the beautiful paper \cite{MS}, except that, as we have
said, we need some new ideas to bypass the obstacles posed by the
sub-Riemannian setting.

A description of the content of the paper is as follows. In Section
\ref{S:IBP} we introduce the relevant geometric setup, and we recall
the main integration by parts theorems from \cite{DGN4} which
constitute the backbone of the paper. In Section \ref{S:growth} we
combine such results with a suitable adaptation of the ideas in
\cite{MS} to establish some general growth results for hypersurfaces
in $\Hn$. A basic new fact is the identity \eqref{remid} in
Proposition \ref{P:remid} which represents the appropriate
sub-Riemannian analogue of \eqref{fi}. Combining it with the
integration by parts we obtain the growth Theorem \ref{T:minsurf},
which concludes Section \ref{S:growth}. Finally, Section
\ref{S:mono} is devoted to proving Theorem \ref{T:monogsintro}.

\section{\textbf{Sub-Riemannian calculus on hypersurfaces}}\label{S:IBP}

In this section we introduce the relevant notation and recall some
basic integration by parts formulas involving the tangential
horizontal gradient on a hypersurface, and the horizontal mean
curvature of the latter, which are special case of some general
formulas discovered in \cite{DGN4}. Such formulas are reminiscent of
the classical one, and in fact they encompass the latter. However,
an important difference is that the ordinary volume form on the
hypersurface $\mS$ is replaced by the $H$-perimeter measure
$d\sigma_H$. Furthermore, they contain additional terms which are
due to the non-trivial commutation relations, which is reflected in
the lack of torsion freeness of the horizontal connection on $\mS$.
Such term prevents the corresponding horizontal Laplace-Beltrami
operator from being formally self-adjoint in $L^2(\mathcal
S,d\sigma_H)$ in general.

We next recall some basic concepts from the sub-Riemannian geometry
of an hypersurface $\mS\subset \Hn$. For a detailed account we refer
the reader to \cite{DGN4}. We consider the Riemannian manifold $M =
\Hn$ with the left-invariant metric tensor with respect to which
$X_1,...,X_{2n}$ is an orthonormal basis, the corresponding
Levi-Civita connection $\nabla$ on $\Hn$, and the horizontal
Levi-Civita connection $\nabla^H$. Let $\Om \subset \Hn$ be a
bounded $C^k$ domain, with $k\geq 2$. We denote by $\n$ the
Riemannian outer unit normal to $\p \Om$, and define the so-called
\emph{angle function} on $\p \Om$ as follows
\begin{equation}\label{ps}
W\ =\ |\bN^H|\ =\ \sqrt{\sum_{j=1}^{2n} <\n, X_j>^2}\ .
\end{equation}

The \emph{characteristic set} of $\Om$, hereafter denoted by
$\Sigma=\Sigma_{\p \Om}$, is the compact subset of $\p \Om$ where
the continuous function $W$ vanishes
\begin{equation}\label{csaf}
\Sigma\ =\ \{p\in \p \Om\mid W(p) = 0\}\ .
\end{equation}

The next definition plays a basic role in sub-Riemannian geometry.

\begin{dfn}\label{D:HGauss}
We define the outer \emph{horizontal normal} on $\p \Om$ as follows
\begin{equation}\label{up1}
\bN^H\ =\ \sum_{j=1}^{2n} <\n,X_j> X_j\  ,
\end{equation}
so that $W = |\bN^H|$. The \emph{horizontal Gauss map} $\nuX$ on $\p
\Om$ is defined by
\begin{equation}\label{up2}
\nuX\ =\ \frac{\bN^H}{|\bN^H|}\ ,\quad\quad\quad\quad \text{on}\quad
\p \Om \setminus \Sigma\ .
\end{equation}

\end{dfn}

Henceforth, we set $<\nuX,X_i> = \pb_i, <\nuX,X_{n+i}> = \qb_i$,
$i=1,...,n$, so that \[ \pb_1^2 + ... + \pb_n^2 + \qb_1^2 + ... +
\qb_n^2 = 1. \] We note that $\Up$ is the projection of the
Riemannian Gauss map on $\p \Om$ onto the horizontal subbundle $H\Hn
\subset T \Hn$. Such projection vanishes only at characteristic
points, and this is why the horizontal Gauss map is not defined on
$\Sigma$. The following definition is taken from \cite{DGN4}, but
the reader should also see \cite{HP2} for a related notion in the
more general setting of vertically rigid spaces.

\begin{dfn}\label{D:HMC}
The \emph{horizontal} or $H$-\emph{mean curvature} of $\p \Om$ at a
point $p_0\in \p \Om \setminus \Sigma$ is defined as
\[
\mathcal H\ =\ \sum_{i=1}^{2n-1} <\nabla^H_{\boldsymbol e_i}
\boldsymbol e_i,\nuX>\ ,
\]
where $\{\boldsymbol e_1,...,\boldsymbol e_{2n-1}\}$ denotes an
orthonormal basis of the horizontal tangent bundle $T_H \p \Om
\overset{def}{=} T \p \Om \cap H\Hn$ on $\p \Om$. If instead $p_0\in
\Sigma$, then we define $\mathcal H(p_0) = \underset{p\to p_0}{\lim}
\mathcal H(p)$, provided that the limit exists and is finite.
\end{dfn}

Given an open set $\Om \subset \Hn$ denote by
\[
\mathcal F(\Om)\ =\ \{\phi = \sum_{j=1}^{2n} \phi_j X_j \in
C^1_0(\Om,H\Hn) \mid ||\phi||_\infty = \underset{ \Om}{\sup}
(\sum_{j=1}^{2n} \phi_j^2)^{1/2} \leq 1\}\ .
\]
Given $\phi = \sum_{j=1}^{2n} \phi_j X_j \in C^1_0(\Om,H\Hn)$, we
let $div_H \phi = \sum_{j=1}^{2n} X_j \phi_j$. The $H$-perimeter of
a measurable set $E\subset \Hn$ with respect to $\Om$ was defined in
\cite{CDG} as
\[
P_H(E;\Om) \ =\ \underset{\phi\in \mathcal F(\Om)}{\sup} \int_{E\cap
\Om} div_H \phi\ dg\ .
\]

If $E$ is a bounded open set of class $C^1$, then the divergence
theorem gives
\[
P_H(E;\Om) = \underset{\phi\in \mathcal F(\Om)}{\sup} \int_{\p E\cap
\Om} \sum_{j=1}^{2n} <\n,X_j> \phi_j\ d\sigma = \underset{\phi\in
\mathcal F(\Om)}{\sup} \int_{\p E\cap \Om} <\bN^H,\phi>  d\sigma =
\int_{\p E\cap \Om} |\bN^H| d\sigma\ ,
\]
where $d\sigma$ is the Riemannian surface measure on $\p E$. It is
clear from this formula that the measure on $\p E$, defined by
\[ \sigma_H(\p E \cap \Om)\ \overset{def}{=}\ P_H(E;\Om)\  \]
on the open sets of $\p E$, is absolutely continuous with respect to
$\sigma$, and its density is represented by the angle function $W$
of $\p E$. We formalize this observation in the following
definition.

\begin{dfn}\label{D:permeas}
Given a bounded domain $E \subset \Hn$ of class $C^1$, with angle
function $W$ as in \eqref{csaf}, we will denote by
\begin{equation}\label{permeasure}
d \sigma_H\ =\ |\Up|\ d\sigma\ =\ W\ d\sigma\ ,
\end{equation}
the $H$-perimeter measure supported on $\p E$.
\end{dfn}

In what follows we will indicate with \[ HT\mS \overset{def}{=} T\mS
\cap H\Hn \] the so-called horizontal tangent bundle of $\mS$.

\begin{dfn}\label{D:horconS}
Let $\mS\subset \Hn$ be a non-characteristic, $C^k$ hypersurface,
$k\geq 2$, then we define the \emph{horizontal connection} on $\mS$
as follows. Let $\nabla^H$ denote the horizontal Levi-Civita
connection in $\Hn$. For every $X,Y\in C^1(\mS;HT\mS)$ we define
\[
\delXY\ =\ \nabla^H_{\oX} \oY\ -\ <\nabla^H_{\oX} \oY,\nuX> \nuX\ ,
\]
where $\oX, \oY$ are any two horizontal vector fields on $\Hn$ such
that $\oX = X$, $\oY = Y$ on $\mS$.
\end{dfn}

One can check that Definition \ref{D:horconS} is well-posed, i.e.,
it is independent of the extensions $\oX, \oY$ of the vector fields
$X, Y$.

\begin{prop}\label{P:nonsymmetry}
For every $X,Y\in C^1(\mS;HT\mS)$ one has
\[
\delXY\ -\ \nabla^{H,\mS}_Y X\ =\ [X,Y]^H\ -\ <[X,Y]^H,\nuX>\nuX\ .
\]
\end{prop}

In the latter identity the notation $[X,Y]^H$ indicates the
projection of the vector field $[X,Y]$ onto the horizontal bundle
$H\Hn$. It is clear from this proposition that the horizontal
connection $\nabla^{H,\mS}$ on $\mS$ is not necessarily torsion
free. This depends on the fact that it is not true in general
that, if $X,Y\in C^1(S;HT\mS)$, then $[X,Y]^H \in C^1(\mS;HT\mS)$.
In the special case of the first Heisenberg group $\HH$ this fact
is true, and we have the following result, see Proposition 7.3 in
\cite{DGN4}.

\begin{prop}\label{P:symmH1}
Given a $C^k$ non-characteristic surface $\mS \subset \HH$, $k\geq
2$, one has $[X,Y]^H\in HT\mS$ for every $X,Y\in C^1(\mS;HT\mS)$,
and therefore the horizontal connection on $\mS$ is torsion free.
\end{prop}

\begin{dfn}\label{D:delta}
Let $\mS$ be as in Definition \ref{D:horconS}. Consider a function
$u \in C^1(\mS)$. We define the \emph{tangential horizontal
gradient} of $u$ as follows
\[
\del u\ \overset{def}{=}\ \nabh \ou\ -\ <\nabh \ou,\nuX>\ \nuX\ ,
\]
where $\ou\in C^1(\bG)$ is such that $\ou = u$ on $\mS$.
\end{dfn}

We are now ready to state the integration by parts formulas from
\cite{DGN4} which constitute the backbone of this paper.

\begin{thrm}[\textbf{Horizontal integration by parts formula}]\label{T:ibp}
Consider a $C^2$ oriented hypersurface $\mathcal S\subset \Hn$. If
$u\in C^1_0(\mS \setminus \Sigma_\mS)$, then we have
\begin{equation}\label{ibp0}
\int_\mathcal S \di u\ d \sigma_H\ =\ \int_\mathcal S u\
\bigg\{\mathcal H\ \nui\ -\ \boldsymbol
c^{H,\mS}_{i}\bigg\}d\sigma_H\ ,\quad\quad\quad i = 1,...,2n\ ,
\end{equation}
where the $C^1$ vector field $\cx = \sum_{i=1}^m \boldsymbol
c^{H,\mS}_{i} X_i$ is given by
\begin{equation}\label{cxHn}
\cx\ =\  \ob\ (\nuX)^\perp\ =\ \ob\ \big(\qb_{1} X_1 + ... + \qb_{n}
X_n - \pb_1 X_{n+1} - ... - \pb_n X_{2n}\big)\ .
\end{equation}
As a consequence, $\cx$ is perpendicular to the horizontal Gauss map
$\nuX$, i.e., one has
\begin{equation}\label{ibp4}
<\cx,\nuX>\ =\ 0\ ,
\end{equation}
and therefore $\cx \in C^1(\mS\setminus \Sigma_\mS, HT\mS)$.
\end{thrm}

\begin{rmrk}\label{R:commutatorHn}
We note explicitly that in view of \eqref{cxHn} we can re-write
\eqref{ibp0} as follows
\begin{equation}\label{ibpHn}
\int_\mathcal S \nabla^{H,\mS} u\ d \sigma_H\ =\ \int_\mathcal S u\
\bigg\{\mathcal H\ \nuX\ -\ \ob\ \nup\bigg\}d\sigma_H\ .
\end{equation}
\end{rmrk}

We have the following notable consequences of Theorem \ref{T:ibp}.

\medskip

\begin{thrm}\label{T:divH}
Let $\mathcal S\subset \Hn$ be a $C^2$ oriented hypersurface, with
characteristic set $\Sigma_\mS$. If $\zeta\in C^1_0(\mS \setminus
\Sigma_\mS,HT\mS)$, then we have
\begin{equation}\label{ibp1}
\int_\mathcal S \left\{div_{H,\mS} \zeta\ +\ <\cx,\zeta>\right\}\ d
\sigma_H\ =\  \int_\mathcal S \mathcal H\ <\zeta,\nuX>\ d\sigma_H\ ,
\end{equation}
where we have let
\[
div_{H,\mS} \zeta\ =\ \sum_{i=1}^{2n} \di \zeta_i\ .
\]
\end{thrm}

We next recall a different integration by parts formula which
involves differentiation along a special combination of the vector
fields $\nuX$ and $T$.

\begin{thrm}[\textbf{Vertical integration by parts formula}]\label{T:T&Y}
Let $\mS\subset \Hn$ be a $C^2$ oriented hypersurface. For every
$f\in C^1(\mS)$,  $g \in C^1_0(\mS\setminus \Sigma_\mS)$, one has
\begin{equation}\label{twoT&Y}
\int_\mS f \ (T  - \ob Y)g\ d\sigma_H\ =\ -\ \int_\mS g \ (T
 - \ob Y)f\ d\sigma_H\ +\ \int_\mS f g \ob\ \mathcal H\
d\sigma_H,
\end{equation}
where we have let $Yf = <\nabla f,\nuX>$.
\end{thrm}


\section{\textbf{Growth formulas for the $H$-perimeter in hypersurfaces}}\label{S:growth}

In this section we present the proof of Theorem \ref{T:monogsintro}.
We begin by introducing some notation. To motivate them we mention
that when we first approached the question of monotonicity of the
$H$-perimeter we asked ourselves whether a result corresponding to
\eqref{vg2} hold for an $H$-minimal hypersurface. More specifically,
in view of the natural $r^{Q-1}$ rescaling of the $H$-perimeter, it
is natural to ask whether for such a hypersurface the function
\[
r\ \to\ \frac{\sigma_H(\mathcal S\cap B(g,r))}{r^{Q-1}}
\]
is monotone non-decreasing.

We next turn our attention to the case of surfaces in $\HH$ which
are in the form of the graphical strips introduced in \cite{DGNP}.
In what follows, we consider functions $\rho \in C^1(\bG)$ and
$\lambda \in C^1(\R)$, to be determined later. We have the following
basic lemma which constitutes a sub-Riemannian counterpart of a
result due to Michael and Simon \cite{MS}.

\begin{lemma}\label{L:ae}
Consider a horizontal vector field $\zeta = \sum_{i=1}^m \zeta_i X_i
\in C^1(\bG,HG)$. Let $\mS\subset \Hn$ be a $C^2$ hypersurface with
empty characteristic locus. Suppose that the level sets of $\rho$
are compact, and let $\lambda$ be non-decreasing, with $\lambda(t)
\equiv 0$ for $t\leq 0$. Given $\psi \in C^1(\mS)$, for every $r>0$
we have
\begin{align}\label{ae0}
& \int_\mS \bigg\{div_{H,\mS} \zeta + <\cx,\zeta> \bigg\} \lambda(r
- \rho)\ \psi\ d\sigma_H\ -\ \int_\mS \lambda'(r - \rho)\ \psi
<\zeta ,\del \rho>\ d\sx
\\
& \leq\ \int_\mS \lambda(r - \rho)\ |\zeta|\ \left\{ |\mathcal H|\
\psi\ +\ |\del \psi|\right\}\ d\sx\ . \notag
\end{align}
In particular, choosing $\psi \equiv 1$ we obtain from \eqref{ae0}
\begin{align}\label{ae00}
& \int_\mS \bigg\{div_{H,\mS}\ \zeta + <\cx,\zeta> \bigg\} \lambda(r
- \rho)\ d\sigma_H\ -\ \int_\mS \lambda'(r - \rho) <\zeta ,\del
\rho>\ d\sx
\\
& \leq\  \int_\mS \lambda(r - \rho)\ |\zeta|\ |\mathcal H|\ d\sx\ .
\notag
\end{align}

\end{lemma}

\begin{proof}[\textbf{Proof}]

For a fixed $r>0$ we define
\begin{equation}\label{ae1}
u\ =\ \zeta_i\ \lambda(r - \rho)\ \psi\ ,
\end{equation}
where $\lambda : \R \to \R$ is non-decreasing, and $\lambda \equiv
0$ for $t\leq 0$. We have
\begin{align}\label{ae2}
\sum_{i=1}^m \di u\ &=\ \left(div_{H,\mS} \zeta\right) \lambda(r -
\rho)\ \psi
\\
& +\ \lambda(r - \rho) <\zeta, \del \psi> \ -\ \lambda'(r - \rho)\
\psi <\zeta, \del \rho>\ . \notag
\end{align}

We now integrate \eqref{ae2} on $\mS$ with respect to the measure
$\sigma_H$. Applying \eqref{ibp0} in Theorem \ref{T:ibp} we obtain
\begin{align}\label{ae4}
& \int_\mS div_{H,S}\ \zeta\ \lambda(r - \rho)\ \psi\ d\sigma_H\ +\
\int_\mS \lambda(r - \rho) <\zeta,\del \psi> d\sx
\\
& -\ \int_\mS \lambda'(r - \rho)\ \psi <\zeta,\del \rho> d\sigma_H\
+\ \int_\mS \lambda(r - \rho)\ \psi <\cx,\zeta> d\sigma_H\
\notag\\
& =\  \int_\mS \mathcal H\ \lambda(r - \rho)\ \psi <\zeta,\nuX> d
\sx\ . \notag
\end{align}

From the identity \eqref{ae4} we easily obtain \eqref{ae0}.

\end{proof}

We next use the formula \eqref{twoT&Y} (with the choice $g(p) =
\lambda(r-\rho(p))$) in Theorem \ref{T:T&Y} to obtain the following
result.

\begin{lemma}\label{L:vibp}
Let $\mS\subset \Hn$ be a $C^2$ hypersurface with empty
characteristic locus. Suppose that the level sets of $\rho$ are
compact, and let $\lambda$ be non-decreasing, with $\lambda(t)
\equiv 0$ for $t\leq 0$. For every $r>0$ we have for any $f\in
C^1(\mS)$
\begin{equation}\label{vibp}
\int_\mS \lambda(r - \rho) \ (T  - \ob Y)f\ d\sigma_H =
  \int_\mS f \ \lambda'(r-\rho) (T
 - \ob Y)\rho\ d\sigma_H +  \int_\mS \lambda(r-\rho) f\ \ob\ \mathcal H\
d\sigma_H\ .
\end{equation}
\end{lemma}

At this point we combine \eqref{ae00} in Lemma \ref{L:ae} with
\eqref{vibp} in Lemma \ref{L:vibp}, obtaining the following basic
result.

\begin{thrm}\label{T:ibpcombined}
Let $\mS\subset \Hn$ be a $C^2$ hypersurface with empty
characteristic locus. Suppose that the level sets of $\rho$ are
compact, and let $\lambda$ be non-decreasing, with $\lambda(t)
\equiv 0$ for $t\leq 0$. For every $r>0$ we have for any $f\in
C^1(\mS)$
\begin{align}\label{ibpc} & \int_\mS
\bigg\{div_{H,\mS}\ \zeta + <\cx,\zeta> + (T  - \ob Y)f\bigg\}
\lambda(r - \rho)\ d\sigma_H
\\
& -\ \int_\mS \lambda'(r - \rho) \bigg\{<\zeta ,\del \rho> + f (T  -
\ob Y) \rho\bigg\}\ d\sx
\notag\\
& \leq\  \int_\mS \lambda(r - \rho)\bigg\{ |\zeta| +  |f|\ob\bigg\}\
|\mathcal H|\ d\sx\ . \notag
\end{align}
In particular, if $\mS$ is $H$-minimal, we obtain from \eqref{ibpc}

\begin{align}\label{ibpcmin}
& \int_\mS \bigg\{div_{H,\mS}\ \zeta + <\cx,\zeta> +  (T  - \ob
Y)f\bigg\} \lambda(r - \rho)\ d\sigma_H
\\
& -\ \int_\mS \lambda'(r - \rho) \bigg\{<\zeta ,\del \rho> + f (T -
\ob Y) \rho\bigg\}\ d\sx
\notag\\
& \leq\  0\ . \notag
\end{align}
\end{thrm}

We now turn to the fundamental question of the choice of the
horizontal vector field $\zeta$ and of the function  $f$ in Theorem
\ref{T:ibpcombined}. With this objective in mind we introduce the
following definition.

\begin{dfn}\label{D:shiftedgendil}
Let $p_0\in \Hn$, then the generator of the non-isotropic dilations
$\{\delta_\lambda\}_{\lambda>0}$ centered at $p_0$ is defined by
\[
Z_{p_0}f(p)\ =\ \sum_{i=1}^n (x_i-x_{0,i})X_i + (y_i-y_{0,i})X_{n+i}
+ [2(t-t_0) + (<x,y_0>-<x_0,y>)] T .
\]
\end{dfn}

Definition \ref{D:shiftedgendil} is motivated by the following
considerations. Let $F\in C^1(\Hn)$, then
\[
Z_{p_0}F(p)\ \overset{def}{=}\ \frac{d}{d\lambda} F(p_0
\delta_\lambda(p_0^{-1} p))\bigg|_{\lambda=1}.
\]
Now
\begin{align*}
p_0 \delta_\lambda(p_0^{-1} p)\ & =\ \bigg(x_0 + \lambda(x-x_0),y
+\lambda(y-y_0),
\\
& t_0  + \lambda^2\left(t-t_0 +\frac{1}{2}(<x,y_0>-<x_0,y>)\right) +
\frac{\lambda}{2}(<x_0,y-y_0> - <y_0,x-x_0>)\bigg)
\end{align*}

A simple calculation now gives
\begin{align}\label{gendil}
& \frac{d}{d\lambda} F(p_0 \delta_\lambda(p_0^{-1}
p))\bigg|_{\lambda=1} \
=\\
& \sum_{i=1}^n (x_i-x_{0,i})\frac{\p F}{\p x_i}(p) +
(y_i-y_{0,i})\frac{\p F}{\p y_i}(p) + [2(t-t_0) +
\frac{1}{2}(<x,y_0>-<x_0,y>)] TF(p) \notag \end{align} If in
\eqref{gendil} we now use the fact that
\[
\frac{\p F}{\p x_i}(p)\ =\ X_i F(p) + \frac{y_i}{2} TF(p),\ \ \
\frac{\p F}{\p y_i}(p)\ =\ X_{n+i} F(p) - \frac{x_i}{2} TF(p),
\]
we easily obtain the formula in Definition \ref{D:shiftedgendil}.

Guided by Definition \ref{D:shiftedgendil}, we now choose the
horizontal vector field $\zeta$ and the function $f$ in Theorem
\ref{T:ibpcombined} as follows
\begin{equation}\label{vfp0}
\zeta(p)\ =\ \sum_{i=1}^n\big((x_i-x_{0,i}) X_i + (y_i - y_{0,i})
X_{n+i}\big), \end{equation}
\begin{equation}\label{f}
f(p)\ =\ 2(t-t_0) + <x,y_0>-<x_0,y>.
\end{equation}

With these choices, we next establish a remarkable identity which
should be considered as the sub-Riemannian counterpart of the above
recalled \eqref{fi}. In what follows, similarly to formula
\eqref{twoT&Y} above, we will use the notation $Yf = <\nabla
f,\nuX>$.

\begin{prop}\label{P:remid}
Fix a point $p_0 = (x_0,y_0,t_0)\in \Hn$ and consider the horizontal
vector field $\zeta\in C^\infty(\Hn,H\Hn)$ given by \eqref{vfp0},
and the function $f\in C^1(\Hn)$ in \eqref{f}, then on any $C^2$
non-characteristic hypersurface $\mS\subset \Hn$ (or on any
hypersurface $\mS$, but away from its characteristic set
$\Sigma_\mS$) one has the identity
\begin{equation}\label{remid}
div_{H,\mS}\ \zeta + <\cx,\zeta> + (T  - \ob Y) f\ \equiv\ Q -1 \ .
\end{equation}
\end{prop}

\begin{proof}[\textbf{Proof}]
We begin by observing that with $\zeta = \sum_{i=1}^n(\zeta_i X_i +
\zeta_{n+i} X_{n+i})$ one has
\[
\di \zeta_i = 1 - \pb_i^2\ ,\ \ \ \nabla^{H,S}_{n+i} \zeta_{n+i} = 1
- \qb_i^2\ ,\ \ i=1,...,n\ .
\]

Therefore,
\begin{equation}\label{Q-3}
div_{H,\mS} \zeta\ =\ \sum_{i=1}^{n} \big(\di \zeta_i +
\nabla^{H,S}_{n+i} \zeta_{n+i}\big)\ =\ 2n - \sum_{i=1}^n(\pb_i^2 +
\qb_i^2)\ \equiv\ 2n - 1\ =\ Q-3\ .
\end{equation}

We now have from \eqref{cxHn}
\[
\cx\ =\ \ob (\qb_1 X_1 + ... + \qb_n X_n - \pb_1 X_{n+1} -... -
\pb_n X_{2n}),
\]
and therefore \[ <\cx,\zeta>\ =\ \ob <z,\nup> - \ob <z_0,\nup>,
\]
where, abusing the notation, we have set $z=\sum_{i=1}^n x_i X_i +
y_i X_{n+i}$, $z_0=\sum_{i=1}^n x_{0,i} X_i + y_{0,i} X_{n+i}$.
 On the other hand,
since $Yt = \frac{1}{2}(x_1 \qb_1 + ... + x_n \qb_n - y_1 \pb_1
-...-y_n \pb_n)$, we have
\[
(T-\ob Y)(2(t-t_0)) = 2 Tt - 2 \ob Yt = 2 - \ob <z,\nup>.
\]
We also have \[ (T-\ob Y)(<x,y_0>-<x_0,y>)\ =\ - \ob
Y(<x,y_0>-<x_0,y>)\ =\ \ob <z_0,\nup>,
\]
and so
\begin{equation}\label{2}
<\cx,\zeta> + (T-\ob Y) f\ \equiv\ 2 \ .
\end{equation}

Combining \eqref{2} with \eqref{Q-3} we obtain \eqref{remid}.
\end{proof}

If we now combine \eqref{ibpcmin} in Theorem \ref{T:ibpcombined}
with Proposition \ref{P:remid}, we obtain the following basic
result.

\begin{thrm}\label{T:minsurf}
Let $\mS\subset \Hn$ be a non-characteristic $H$-minimal surface,
then with $\zeta$ as in \eqref{vfp0} and $f$ as in \eqref{f}, one
has for any $p_0=(x_0,y_0,t_0)\in \Hn$

\begin{align}\label{ibpcmin2}
& (Q-1) \int_\mS \lambda(r - \rho)\ d\sigma_H
\\
& -\ \int_\mS \lambda'(r - \rho) \bigg\{<\zeta ,\del \rho> + f (T -
\ob Y) \rho\bigg\}\ d\sx
\notag\\
& \leq\  0\ . \notag
\end{align}
\end{thrm}

\section{\textbf{Monotonicity for graphical strips}}\label{S:mono}

In this section we obtain an interesting consequence of Theorem
\ref{T:minsurf} by proving an intrinsic monotonicity property
similar to \eqref{vg2} for a remarkable class of $H$-minimal
surfaces in the Heisenberg group $\HH$. Such surfaces, called
\emph{graphical strips} in \cite{DGNP}, have been introduced in
connection with the solution of the sub-Riemannian Bernstein problem
in \cite{DGNP}. The following result is part of Theorem 1.5 in
\cite{DGNP}.

\begin{prop}\label{P:gs}
Every $C^2$ graphical strip is an $H$-minimal surface with empty
characteristic locus. \end{prop}

We recall that, given a $C^1$ surface $\mS\subset \HH$, the
characteristic locus of $\mS$, henceforth denoted by $\Sigma_\mS$,
is the collection of all points $p\in \mS$ at which $H_p\HH = T_p
\mS$, where $H_p\HH$ denotes the fiber at $p$ of the horizontal
bundle of $\HH$. One fundamental aspect of graphical strips is
represented by the following result, which constitutes one of the
two central results in \cite{DGNP}. In order to state it we mention
that $\nuX$ indicates the horizontal Gauss map of $\mS$, which is
well defined away from the characteristic locus $\Sigma_\mS$ of
$\mS$. By $\mathcal V^H_{II}(\mS; \mathcal X)$ we denote the second
variation of the $H$-perimeter with respect to a deformation of
$\mS$ in the direction of the vector field $\mathcal X$. An
$H$-minimal surface $\mS$ with empty characteristic locus is called
\emph{stable} if $\sv \geq 0$ for every compactly supported
$\mathcal X = a X_1 + b X_2 + k T$. Otherwise, it is called
unstable. We note that, since thanks to Proposition \ref{P:gs} every
graphical strip has empty characteristic locus, the horizontal Gauss
map $\nuX$ of such a surface is globally defined.

\begin{thrm}\label{T:perminI}
Let $\mS$ be a $C^2$ strict graphical strip, then $\mS$ is unstable.
In fact, there exists a continuum of $h\in C^2_0(\mS)$, for which
$\mathcal V^H_{II}(\mS; h \nuX) <  0$.
\end{thrm}

The following theorem constitutes the second main result in
\cite{DGNP}. It underscores the central relevance of graphical
strips in the study of $H$-minimal surfaces in $\HH$.

\begin{thrm}\label{T:reductionI} Let
$\mS\subset \HH$ be an $H$-minimal entire graph of class $C^2$, with
empty characteristic locus, and that is not itself a vertical plane
\begin{equation}\label{vp}
 \Pi_0  = \{(x,y,t)\in \HH\mid ax + by = \gamma_0\},
\end{equation}
then there exists a strict graphical strip $\mS_0 \subset \mS$.
\end{thrm}

By combining Theorems \ref{T:perminI} and \ref{T:reductionI} the
following solution of the sub-Riemannian Bernstein problem was
obtained in \cite{DGNP}.

\begin{thrm}[\textbf{of Bernstein type}]\label{T:instability}
In $\HH$ the only $C^2$ stable $H$-minimal entire graphs, with empty
characteristic locus, are the vertical planes \eqref{vp}.
\end{thrm}

In connection with the stability assumption in Theorem
\ref{T:instability} it should be emphasized that, without it, the
theorem is false. This central aspect of the problem was first
discovered in \cite{DGN5} where it was shown that the non-planar
$H$-minimal surface $\mS = \{(x,y,t)\in \HH\mid x = yt\}$ (which is
easily seen to be an entire strict graphical strip) is unstable.

Henceforth, given a $C^1$ surface $\mS\subset \HH$ we will indicate
with $\sigma_H$ the horizontal perimeter measure on $\mS$. We
emphasize (see for instance \cite{DGN3}), that such measure scales
according to the following equation

\begin{equation}\label{scalingsE}
\sigma_H(\delta_\lambda(E))\ =\ \lambda^{Q-1} \sigma_H(E),
\end{equation}
with respect to the non-isotropic group dilations
$\delta_\lambda(x,y,t) = (\lambda x,\lambda y, \lambda^2 t)$. Here,
the number $Q = 2n+2$ represents the homogeneous dimension of $\Hn$
associated with the dilations $\{\delta_\lambda\}_{\lambda>0}$. For
instance, when $n=1$, then we have $Q = 4$.

The main result of the present section is the following theorem.

\begin{thrm}\label{T:monogs0}
Let $\mS\subset \HH$ be a $C^2$ graphical strip, then for every
$p_0=(0,0,t_0)\in \mS$ the function
\[
r\ \to\ \frac{\sigma_H(\mS\cap B(p_0,r))}{r^{Q-1}},\ \ \ \ r>0,
\]
is monotone non-decreasing. Moreover, there exists $\omega>0$ such
that
\[
\sigma_H(\mS\cap B(p_0,r)) \geq \omega r^{Q-1},\ \ \ \text{for
every}\ r>0.
\]
\end{thrm}

In the statement of Theorem \ref{T:monogs0} we have denoted by
$B(p_0,r) = \{p\in \Hn \mid d(p,p_0)<r\}$, where $d(p,p_0) =
N(p_0^{-1} p)$ represents the gauge distance on $\Hn$ defined via
the Koranyi-Folland gauge function $N(p) = (|z|^4 + 16 t^2)^{1/4}$,
$p = (z,t)\in \Hn$.

We now specialize the choice of the function $\rho$ in Theorem
\ref{T:minsurf} by letting $\rho(p) = N(p_0^{-1} p)$. Of course,
this is not the only possible choice of $\rho$, but at the moment we
will not further investigate this question since we plan to return
to it in a future study.

Notice that we can write
\begin{equation}\label{ag}
\rho(p) = \left[\big((x-x_0)^2 + (y-y_0)^2\big)^2 + 4 \big(2(t-t_0)
+ (xy_0 - x_0y)\big)^2\right]^{1/4}\ .
\end{equation}

A simple calculation gives \begin{align}\label{gg}  X_1\rho & =
\rho^{-3}\bigg[(x-x_0) |z-z_0|^2 - 2(y-y_0)\big(2(t-t_0) + (xy_0 -
x_0y)\big)\bigg], \\ \label{gg2}
 X_2\rho & =
\rho^{-3}\bigg[(y-y_0) |z-z_0|^2 + 2(x-x_0)\big(2(t-t_0) + (xy_0 -
x_0y)\big)\bigg],
\\ \label{gg3}
T\rho & = \rho^{-3} 4 \bigg[2(t-t_0) + (xy_0 - x_0y)\bigg].
\end{align}

From \eqref{gg}, \eqref{gg2}, \eqref{gg3} we obtain with $\zeta$ and
$f$ as in \eqref{vfp0}, \eqref{f} respectively,
\begin{equation}\label{sg}
<\zeta,\nabh \rho> + f T\rho\ =\ \rho.
\end{equation}

On the other hand, we have from the expression of the horizontal
covariant derivative on $\mS$
\[ <\zeta ,\del \rho>\ =\ <\zeta,\nabh
\rho> - <\nabh \rho,\nuX><\zeta,\nuX>.
\]

Using \eqref{sg} we find
\begin{align}\label{uffa}
& <\zeta ,\del \rho> + f (T - \ob Y) \rho
\\
& =\  <\zeta,\nabh \rho> + f T\rho \notag\\
& - <\nabh \rho,\nuX><\zeta,\nuX> - \ob f
Y\rho \notag\\
& =\ \rho - <\nabh \rho,\nuX> \bigg(<\zeta,\nuX> + \ob  f \bigg),
\notag
\end{align}
where in the last equality we have used the fact that $Y\rho =
<\nabh \rho,\nuX>$.

The next result provides a fundamental estimate. It is at this point
that we use the special structural assumption that $\mS$ be a
graphical strip in $\HH$.

\begin{lemma}\label{L:crucial}
Let $\mS\subset \HH$ be a $C^2$ graphical strip. Let $p_0 =
(0,0,t_0)\in \mS$, then with $\zeta$ as in \eqref{vfp0} and $f$ as
in \eqref{f}, one has
\[
\underset{\mS\cap B(p_0,r)}{\sup}\ \left|<\zeta ,\del \rho> + f (T -
\ob Y) \rho\right|\ \leq\ r.
\]
\end{lemma}

\begin{proof}[\textbf{Proof}]
In view of \eqref{uffa}, proving the lemma is equivalent to showing
\[
\underset{\mS\cap B(p_0,r)}{\sup}\ \left|\rho - <\nabh \rho,\nuX>
\bigg(<\zeta,\nuX> + \ob  f \bigg)\right|\ \leq\ r.
\]

Without loss of generality we assume that
\[
\mS = \{(x,y,t)\in \HH\mid (y,t)\in \R\times I, x = y G(t)\},
\]
for some $G\in C^2(I)$, such that $G'(t)\geq 0$ for every $t\in I$.
We next recall some calculations from \cite{DGNP}. It is obvious
from the definition that $\mS$ is a $C^2$ graph over the
$(y,t)$-plane.  We can use the global defining function
\begin{equation}\label{df}
\phi(x,y,t)\ =\ x - y G(t)\ , \end{equation} and assume that $\mS$
is oriented in such a way that a non-unit Riemannian normal on $\mS$
be given by $\bN = \nabla \phi = (X_1\phi)X_1 + (X_2\phi) X_2 +
(T\phi) T$. We thus find
\begin{equation}\label{pqt}
p = X_1\phi = 1 + \frac{y^2}{2} G'(t)\ ,\quad \quad q = X_2 \phi = -
G(t) - \frac{xy}{2} G'(t)\ ,\quad\quad \omega\ = T\phi = - y G'(t)\
.
\end{equation}
Since $p\geq 1>0$, we see from \eqref{pqt} that $\Sigma_\mS =
\varnothing$.

From now on, to simplify the notation, we will omit the variable $t$
in all expressions involving $G(t), G'(t)$. The second equation in
\eqref{pqt} becomes on $\mS$
\begin{equation}\label{qons}
q\ =\ -\ G\ \left(1\ +\ \frac{y^2}{2} G'\right)\ .
\end{equation}

We thus find on $\mS$
\begin{equation}\label{Wons}
W\ =\ \sqrt{p^2 + q^2}\ =\ \sqrt{1 + G^2} \left(1\ +\ \frac{y^2}{2}
G'\right)\ .
\end{equation}

The equations \eqref{pqt}, \eqref{qons} and \eqref{Wons} give on
$\mS$
\begin{equation}\label{pbqbS}
\pb = \frac{1}{\sqrt{1 + G^2}}\ ,\ \ \qb =  - \frac{G}{\sqrt{1 +
G^2}}\ ,\ \ \ \ob =  - \frac{y G'}{\sqrt{1 + G^2} \left(1\ +\
\frac{y^2}{2} G'\right)}.
\end{equation}

We thus have on $\mS$
\begin{equation}\label{xypq}
x \qb - y \pb = - \left\{\frac{y G^2}{\sqrt{1 + G^2}} +
\frac{y}{\sqrt{1 + G^2}}\right\} = - y \sqrt{1 + G^2}\ ,
\end{equation}
and also
\begin{equation}\label{zero}
x\pb + y \qb \ =\ \frac{y\,G(t)}{\sqrt{1 + G(t)^2}} -
\frac{y\,G(t)}{\sqrt{1 + G(t)^2}} \ =\ 0\ .
\end{equation}

On the other hand, if $p_0 = (x_0,y_0,t_0)\in \mS$, we must have
$x_0 = y_0 G(t_0)$, and therefore
\begin{equation}\label{xypq0}
x_0 \qb - y_0 \pb = - y_0\left\{\frac{G(t_0)G}{\sqrt{1 + G^2}} +
\frac{1}{\sqrt{1 + G^2}}\right\} = - y_0 \frac{1 + G(t_0)G}{\sqrt{1
+ G^2}},
\end{equation}
and also
\begin{equation}\label{finalest0}
x_0\pb + y_0 \qb \ =\ - y_0 \frac{G - G(t_0)}{\sqrt{1 + G^2}}\ .
\end{equation}

We also have on $\mS$
\begin{equation}\label{hp}
xy_0 - x_0y\ =\ y_0 y(G - G(t_0))\ .
\end{equation}

Combining \eqref{zero} and \eqref{finalest0} we find
\begin{equation}\label{fe3}
<\zeta,\nuX>\ =\ (x-x_0)\pb + (y - y_0) \qb\ =\ y_0 \frac{G -
G(t_0)}{\sqrt{1 + G^2}}\ .
\end{equation}

From \eqref{pbqbS}, \eqref{hp} we have \begin{align}\label{ob} \ob
(2(t-t_0) + (xy_0 - x_0y))& \ =\ - 2 y_0 \frac{\frac{y^2}{2} G'(G -
G(t_0))}{\sqrt{1 + G^2} \left(1\ +\ \frac{y^2}{2} G'\right)}
 \\
 & - \frac{2 (t-t_0)y G'}{\sqrt{1 + G^2} \left(1\ +\
\frac{y^2}{2} G'\right)}\ \notag
\end{align}

Combining \eqref{fe3} and \eqref{ob} we find
\begin{align}\label{fe4}
& <\zeta,\nuX> + f \ob  = y_0 \frac{G - G(t_0)}{\sqrt{1 + G^2}}
\\
&  - 2 y_0 \frac{\frac{y^2}{2} G'(G - G(t_0))}{\sqrt{1 + G^2}
\left(1\ +\ \frac{y^2}{2} G'\right)} - \frac{2 (t-t_0)y G'}{\sqrt{1
+ G^2} \left(1\ +\ \frac{y^2}{2} G'\right)}. \notag
\end{align}

When $x_0=y_0=0$, and therefore $p_0 = (0,0,t_0)$, we obtain from
\eqref{fe4}
\begin{equation}\label{fe5}
<\zeta,\nuX> + f \ob  =  - \frac{2 (t-t_0)y G'}{\sqrt{1 + G^2}
\left(1\ +\ \frac{y^2}{2} G'\right)}.
\end{equation}

Next, we observe that we have on $\mS$
\[
|z|^2 = y^2(1 + G^2),\ \ x|z|^2 = y^3G(1+G^2),\ \ y|z|^2 = y^3(1 +
G^2).
\]

If we use these formulas in \eqref{gg}, \eqref{gg2}, in combination
with \eqref{pbqbS}, we obtain
\begin{equation}\label{fe6}
<\nabh \rho,\nuX> = - \frac{4 y(t-t_0) (1 + G^2)}{\rho^3
\sqrt{1+G^2}}.
\end{equation}

Combining equations \eqref{fe5}, \eqref{fe6} we find
\begin{equation}\label{fe7}
<\nabh \rho,\nuX>\big(<\zeta,\nuX> + f \ob\big) = \frac{16 (t-t_0)^2
\frac{y^2}{2} G'}{\rho^3 \left(1\ +\ \frac{y^2}{2} G'\right)}.
\end{equation}

Since on $\mS$ we have
\[
\rho^4 = (x^2 + y^2)^2 +16 (t-t_0)^2 = y^4(1 + G^2)^2 + 16
(t-t_0)^2,
\]
from this equation and from \eqref{fe7} it is at this point easy to
check that on $\mS$ one has
\[
\rho - <\nabh \rho,\nuX> \bigg(<\zeta,\nuX> + \ob  f \bigg) \geq 0.
\]

Since from \eqref{fe7} again we see that $<\nabh
\rho,\nuX>\big(<\zeta,\nuX> + f \ob\big)\geq 0$, we finally obtain
\[
\left|\rho - <\nabh \rho,\nuX> \bigg(<\zeta,\nuX> + \ob  f
\bigg)\right| = \rho - <\nabh \rho,\nuX> \bigg(<\zeta,\nuX> + \ob  f
\bigg)\leq \rho,
\]
which, in particular, proves the lemma.
\end{proof}

We can now prove the main result in this section.

\begin{proof}[\textbf{Proof of Theorem \ref{T:monogs0}}]
We define
\begin{equation}\label{P}
\mathcal P(r) = \int_S \lambda(r - \rho) d\sigma_H.
\end{equation}
We easily find
\[
\frac{d}{dr} \left(\frac{\mathcal P(r)}{r^{Q -1}}\right) =
\frac{1}{r^Q} \bigg(r \mathcal P'(r) - (Q-1) \mathcal P(r)\bigg)\ .
\]
We next recall that for any $p_0=(x_0,y_0,t_0)\in \Hn$ one has from
\eqref{ibpcmin2},

\begin{align}\label{ibpcmin3}
(Q-1) \mathcal P(r) -\ \int_\mS \lambda'(r - \rho) \bigg\{<\zeta
,\del \rho> + f (T - \ob Y) \rho\bigg\}\ d\sx \leq  0,
\end{align}
where $\zeta$ is as in \eqref{vfp0} and $f$ as in \eqref{f}.

At this point the crucial Lemma \ref{L:crucial} enters the picture.
In it we have proved that on the set $B(p_0,r) = \{\rho<r\}$ one has
\begin{equation}\label{crucial}
\left|<\zeta ,\del \rho> + f (T  - \ob Y) \rho\right|\ \leq\ r\ .
\end{equation}

Then from \eqref{crucial}, the fact that $\lambda'(r -\rho)\geq 0$
and from \eqref{ibpcmin3} we can conclude that
\[
\frac{d}{dr} \left(\frac{\mathcal P(r)}{r^{Q -1}}\right) =
\frac{1}{r^Q} \bigg(r \mathcal P'(r) - (Q-1) \mathcal P(r)\bigg)\
\geq\ 0\ .
\]

We now fix $0<r_1<r_2<\infty$ and integrate the latter inequality on
the interval $(r_1,r_2)$ obtaining
\begin{align}\label{fe8}
0 & \leq \int_{r_1}^{r_2} \frac{d}{dr} \left(\frac{\mathcal
P(r)}{r^{Q -1}}\right) dr = \frac{\mathcal P(r_2)}{r_2^{Q -1}} -
\frac{\mathcal P(r_1)}{r_1^{Q -1}}
\\
& = \frac{1}{r_2^{Q -1}} \int_\mS \lambda(r_2 - \rho) d\sigma_H -
\frac{1}{r_1^{Q -1}} \int_\mS \lambda(r_1 - \rho) d\sigma_H \notag
\end{align}

At this point we fix arbitrarily $0<\epsilon <r_1$, and choose a
non-decreasing $0\leq \lambda(s)\leq 1$, with $\lambda\equiv 0$ if
$s\leq 0$, $\lambda \equiv 1$ if $s\geq \epsilon$. With this choice
we obtain from \eqref{fe8}
\begin{align}\label{fe9}
0 & \leq \frac{1}{r_2^{Q -1}} \int_{\mS\cap B(p_0,r_2)} \lambda(r_2
- \rho) d\sigma_H - \frac{1}{r_1^{Q -1}} \int_{\mS\cap
B(p_0,r_1-\epsilon)} \lambda(r_1 - \rho) d\sigma_H \\
& - \frac{1}{r_1^{Q -1}} \int_{\mS\cap [B(p_0,r_1)\setminus
B(p_0,r_1-\epsilon)]} \lambda(r_1 - \rho) d\sigma_H \notag\\
& \leq \frac{\sigma_H(\mS\cap B(p_0,r_2))}{r_2^{Q -1}} -
\frac{\sigma_H(\mS\cap B(p_0,r_1 - \epsilon))}{r_1^{Q -1}}.\notag
\end{align}

Letting $\epsilon \to 0$ we reach the conclusion
\[
\frac{\sigma_H(\mS\cap B(p_0,r_1))}{r_1^{Q -1}} \leq
\frac{\sigma_H(\mS\cap B(p_0,r_2))}{r_2^{Q -1}}.
\]
\end{proof}

According to Theorem \ref{T:monogs0} the limit
\[
\underset{r\to 0^+}{\lim}\ \frac{\sigma_H(\mS \cap
B(p_0,r))}{r^{Q-1}}
\]
exists. In the next proposition we show that such limit is actually
positive.

\begin{prop}\label{P:atzero}
Let $\mathcal S$ be a graphical strip, that is,
\[
\mathcal S\ =\ \{(x,y,t)\,|\, x = y\,G(t)\}\qquad\text{where } G \in
C^1(\mathbb R),\quad G'(t) \geq 0 \text{ for all } t \in \mathbb R\
,
\]
then for every $p_0 = (0,0,t_0)\in \mS$ we have
\[
\underset{r \to 0^+}{\lim}\, \frac{\sigma_H(\mathcal S \cap
B(p_0,r))}{r^3} \ =\ \int_0^1 (1 - \tau^2)^\frac{1}{4}\,d\tau\ > 0 \
.
\]
Note that this limit is independent of $G(t)$.
\end{prop}

\begin{proof}[\textbf{Proof}]
Let $\phi$ be as in \eqref{df}. We then have
\[
\mathcal S \cap B(p_0,r) \ =\ \{(x,y,t)\in \HH\,|\, x = y\,G(t)\ ,\
y^4(1 + G(t)^2)^2 + 16\,(t-t_0)^2 < r^4 \}\ .
\]

\begin{align*}
& |X\phi|\ =\ \left(1 + \frac{y^2}{2}G'(t)\right)\,\sqrt{1 + G(t)^2}\ , \\
&
\end{align*}
Hence
\begin{align}\label{quot}
& \frac{\sigma_H(\mathcal S \cap B(p_0,r))}{r^3}\ =\ \frac{1}{r^3}\,
\int_{\mathcal S \cap B(p_0,r)}\frac{|X\phi|}{|\nabla \phi|}\,d\sigma \\
\notag &\ =\ \frac{1}{r^3}\,\int_{\{(y,t)\,|\,y^4(1+G(t)^2)^2 + 16
(t-t_0)^2 < r^4\}} \left(1 + \frac{y^2}{2}G'(t)\right)\,\sqrt{1 +
G(t)^2}\ dy\,dt
\\
\notag &\ =\ \frac{1}{r^3}\,\int_{t_0-\frac{r^2}{4}}^{t_0 +
\frac{r^2}{4}} \sqrt{1 + G(t)^2} \left(\int_{-\frac{(r^4 -
16(t-t_0)^2)^\frac{1}{4}}{\sqrt{1 + G(t)^2}}}^{\frac{(r^4 -
16(t-t_0)^2)^\frac{1}{4}}{\sqrt{1 + G(t)^2}}}
\ \left(1 + \frac{y^2}{2}G'(t)\right)\,dy\right)\,dt \\
\notag &\ =\
\frac{2}{r^3}\,\int_{t_0-\frac{r^2}{4}}^{t_0+\frac{r^2}{4}} \sqrt{1
+ G(t)^2}
\left\{ \frac{(r^4 - 16(t-t_0)^2)^\frac{1}{4}}{\sqrt{1 + G(t)^2}}\ +\ \frac{G'(t)}{6}\frac{(r^4 - 16(t-t_0)^2)^\frac{3}{4}}{(1 + G(t)^2)^\frac{3}{2}}\right\}\,dt \\
\notag &\ =\ \frac{2}{r^3}\,\int_{t_0
-\frac{r^2}{4}}^{t_0+\frac{r^2}{4}} (r^4 - 16(t-t_0)^2)^\frac{1}{4}
\ +\
\frac{G'(t)}{6}\frac{(r^4 - 16(t-t_0)^2)^\frac{3}{4}}{1 + G(t)^2}\,dt \\
\notag &\ =\ \frac{2}{r^3}\,\int_{t_0
-\frac{r^2}{4}}^{t_0+\frac{r^2}{4}} (r^4 - 16(t-t_0)^2)^\frac{1}{4}\
dt \ +\ \frac{2}{r^3}\,\int_{t_0-\frac{r^2}{4}}^{t_0+\frac{r^2}{4}}
\frac{G'(t)}{6(1 + G(t)^2)} \,(r^4 - 16(t-t_0)^2)^\frac{3}{4}\ dt\ .
\notag
\end{align}
To continue we make the change of variable $t -t_0=
\frac{r^2}{4}\tau$ and analyze the following two terms.

\begin{equation}\label{T1}
\frac{2}{r^3}\,\int_{t_0-\frac{r^2}{4}}^{t_0+\frac{r^2}{4}} (r^4 -
16(t-t_0)^2)^\frac{1}{4}\ dt \ =\ \frac{2}{r^3}\int_{-1}^1 r\,(1 -
\tau^2)^\frac{1}{4}\ \frac{r^2}{4}\ d\tau \ =\ \int_0^1 (1 -
\tau^2)^\frac{1}{4}\ d\tau
\end{equation}

\begin{align}\label{T2}
& \frac{2}{r^3}\,\int_{t_0-\frac{r^2}{4}}^{t_0+\frac{r^2}{4}}
\frac{G'(t)}{6(1 + G(t)^2)}
\,(r^4 - 16(t-t_0)^2)^\frac{3}{4}\ dt \\
\notag & \qquad \qquad \ =\
\frac{1}{3r^3}\int_{-1}^1 \frac{G'(t_0+ r^2\tau/4)}{1 + G(t_0+r^2\,\tau/4)^2}\ r^3\,(1 - \tau^2)^\frac{3}{4}\ \frac{r^2}{4}\ d\tau \\
\notag & \qquad \qquad \ =\ \frac{r^2}{12} \int_{-1}^1
\frac{G'(t_0+r^2\tau/4)}{1 + G(t_0+r^2\,\tau/4)^2}\ (1 -
\tau^2)^\frac{3}{4}\ d\tau \ . \notag
\end{align}
Using \eqref{T1} and \eqref{T2} in \eqref{quot} and Lebesgue
dominated convergence theorem, we obtain

\begin{align}\label{limit0}
& \underset{r \to 0^+}{\lim}\,
\frac{\sigma_H(\mathcal S \cap B(p_0,r))}{r^3} \\
\notag & \qquad\qquad \ =\ \underset{r \to 0^+}{\lim}\, \int_0^1 (1
- \tau^2)^\frac{1}{4}\,d\tau \ +\ \underset{r \to 0^+}{\lim}\,
 \frac{r^2}{12} \int_{-1}^1 \frac{G'(t_0+r^2\tau/4)}{1 + G(t_0+r^2\,\tau/4)^2}\ (1 - \tau^2)^\frac{3}{4}\ d\tau \\
\notag & \qquad\qquad \ =\ \int_0^1 (1 - \tau^2)^\frac{1}{4}\,d\tau
\ +\ \left(\underset{r \to 0^+}{\lim}\,  \frac{r^2}{12} \right)
\int_{-1}^1 (1 - \tau^2)^\frac{3}{4}\ \underset{r \to 0^+}{\lim}\
\frac{G'(t_0+r^2\tau/4)}{1 + G(t_0+r^2\,\tau/4)^2}
\ d\tau \\
\notag & \qquad\qquad \ =\ \int_0^1 (1 - \tau^2)^\frac{1}{4}\,d\tau\
. \notag
\end{align}
\end{proof}

Using \eqref{quot}, \eqref{T1} and \eqref{T2} we can also compute
and obtain

\begin{align}\label{limit1}
& \underset{r \to \infty}{\lim}\,
\frac{\sigma_H(\mathcal S \cap B(p_0,r))}{r^3} \\
\notag & \qquad\qquad \ =\ \int_0^1 (1 - \tau^2)^\frac{1}{4}\,d\tau
\ +\ \underset{r \to \infty}{\lim}\ \int_{-1}^1
\frac{r^2\,G'(t_0+r^2\tau/4)}{12(1 + G(t_0+r^2\,\tau/4)^2)}\ (1 -
\tau^2)^\frac{3}{4}\ \ d\tau \ .
\end{align}

Of course, the above limit may or may not be finite.

At this point, combining Theorem \ref{T:monogs0} and Proposition
\ref{P:atzero} we obtain the maximum sub-Riemannian volume growth of
graphical strips at infinity.

\begin{cor}\label{C:mvg}
Let $\mS\subset \HH$ be a graphical strip, then for every
$p_0=(0,0,t_0)\in \mS$, and every $r>0$ one has
\[
\sigma_H(\mS\cap B(p_0,r)) \geq \omega r^{Q-1},
\]
where we have set $\omega = \int_0^1 (1 -
\tau^2)^\frac{1}{4}\,d\tau$.
\end{cor}

\end{document}